\documentclass[11pt,a4paper]{article}

\pdfoutput=1
\usepackage[utf8]{inputenc}
\usepackage{amsmath,bm}
\usepackage{amsfonts}
\usepackage{amssymb}
\usepackage{upgreek}
\usepackage{fontenc}
\usepackage{comment}
\usepackage{graphicx}
\usepackage{accents}
\usepackage{trimclip}
\usepackage[title]{appendix}
\usepackage{caption}
\usepackage{subcaption} 

\newcommand{\myPiS}[0]{\ds \frac{\pi}{6}}

\newcommand {\der} [2] {\ds {\frac{\partial #1}{\partial #2}}}

\newcommand {\ds}{\displaystyle}

\let \oldnabla \nabla
\renewcommand {\nabla} [0] {\bm{\oldnabla}}

\newcommand{\rn}[1]{%
  \textup{\uppercase\expandafter{\romannumeral#1}}%
}

\usepackage{stackengine}
\stackMath
\newcommand\tenq[2][1]{%
 \def\useanchorwidth{T}%
  \ifnum#1>1%
    \stackunder[0pt]{\tenq[\numexpr#1-1\relax]{#2}}{\scriptscriptstyle\sim}%
  \else%
    \stackunder[1pt]{#2}{\scriptscriptstyle\sim}%
  \fi%
}
\stackMath
\newcommand\mybar[2][1]{%
 \def\useanchorwidth{T}%
  \ifnum#1>1%
    \stackunder[0pt]{\mybar[\numexpr#1-1\relax]{#2}}{\scriptscriptstyle\textminus}%
  \else%
    \stackunder[1pt]{#2}{\scriptscriptstyle\textminus}%
  \fi%
}

\newcommand{\tensf}[1] 	{\underaccent \bar{{\bm{#1}}}}
\newcommand{\tenss}[1] 	{\underaccent \tilde{{\bm{#1}}}}
\newcommand{\tenst}[1] 	{\underaccent \bar{\underaccent \tilde{{\bm{#1}}}}}
\newcommand{\tensff}[1] 	{\underaccent \tilde{\underaccent \tilde{{\bm{#1}}}}}
\newcommand{\bv}[1] 	{\underaccent \bar{\bm{#1}}}
\newcommand{\bvh}[1] 	{\underaccent \bar{\hat{\bm{#1}}}}

\newcommand{\cdotsd} {\cdot \! \cdot \,}

\begin{document}
\title{A micropolar isotropic plasticity formulation for non-associated flow rule  and softening featuring multiple classical yield criteria   \newline Part I - Theory}
\author{Andrea Panteghini, Rocco Lagioia}
\maketitle

\section{Abstract}
The Cosserat continuum is used in this paper to regularize the ill-posed governing equations of the Cauchy/Maxwell continuum.
Most available constitutive models adopt yield and plastic potential surfaces with a circular deviatoric section. This is  a too crude an approximation which hinders the application of  the Cosserat continuum into practice, particularly in the geotechnical domain.
An elasto-plastic constitutive model for the linear formulation of the  Cosserat continuum  is here presented, which features non-associated flow and hardening/softening behaviour, whilst linear hyper-elasticity is adopted to reproduce the recoverable response.
For the formulation  of the yield and plastic potential functions,  a definition of the  \textit{equivalent von Mises stress} is used which is based on Hencky's interpretation of the von Mises criterion and also on the theory of representations. 
The dependency on the Lode's angle of both the yield and plastic potential functions is introduced through the adoption of a recently proposed \textit{Generalized classical} criterion, which rigorously defines most of the classical yield and failure criteria.

\section{Introduction}
Elasto-plastic constitutive models  of the classical Cauchy/Maxwell continuum are widely used to analyse  structural and geotechnical problems.  
The need of accounting for important features of the material behaviour usually requires, particularly for soils, the introduction of non-associativeness and/or softening in the constitutive formulation,  which unfortunately results in serious consequences.

It has long been known that numerical analyses of boundary value problems exhibit non uniqueness of the solution which results in mesh sensitivity when the material is characterized by softening behaviour.  This manifests itself as mesh dependent results,  which do not converge to the correct solution when the domain discretisation  is refined.  Moreover   the  width of the  band with localized deformations depends on the size of the elements (see e.g.  De Borst \cite{DeBorst1991},  De Borst and Sluys \cite{DeBorst1991a},  Stefanou et al.  \cite{Stefanou2019}).  
The adoption of a non-associated flow rule too results in mesh dependency and also in structural softening even if the material exhibits perfect plastic response (e.g.  Sabet and De Borst \cite{Sabet2019}).

However one of the most severe practical consequences for the analyst is that numerical analyses of boundary value problems are very much prone to crushing and indeed there are a number of problems for which the analyses do not perform even a single step when non-associativity and/or softening is introduced.
It is now well established  that mesh sensitivity is the manifestation in numerical methods of the loss of well-posedness of the governing equations  which for static problems  loose ellipticity when localisations occur (e.g.  Sabet and de Borst  \cite{Sabet2019}).

The introduction of a characteristic length is a valid approach for regularizing the problem, and this is indeed what a Cosserat continuum does via the introduction of additional rotational degrees of freedom.

Most available constitutive models for the  Cosserat continuum adopt yield and plastic potential surfaces  defined by classical yield/failure criteria with a circular deviatoric section, such as the von Mises and the Drucker-Prager ( e.g.  Altenbach and Eremeyev \cite{Altenbach2014},  Russo et al. \cite{Russo2020},  Sabet and De Borst \cite{Sabet2019}). 
For many applications, such as most of those in the   geotechnical  domain,  that is too much of  a crude  assumption,  which, depending on how parameters are calibrated,  results in excessively  unconservative or overconservative predictions.  This considerably hinders the application of the Cosserat continuum into practical engineering problems.

An elasto-plastic constitutive model for the Cosserat continuum is here presented featuring linear hyper-elastic behaviour,  non-associated flow  and hardening/softening response. 
The yield and the plastic potential surfaces are those of the main classical yield and failure criteria for metals and soils which are selected by choosing the appropriate set of parameters. This is achieved by adopting the mathematical formulation of those criteria recently presented by Lagioia and Panteghini \cite{LP2016}.
The dependence of the yield and plastic potential functions on the Lode's angle is hence introduced in the model. 
Moreover given the generality of the proposed model,  other classical constitutive relatioships such as the Modified Cam-Clay  (e.g. Roscoe and Burland \cite{Roscoe1968},  Panteghini and Lagioia \cite{Panteghini2018})  can be easily introduced.

\section{Notation}
In this manuscript we will follow the classical continuum mechanics sign convention, whereby tensile stresses and strains are positive.

Tensors of the first, second, third and fourth orders  will be referred to in the compact representation using bold letters of any case,  but for the purpose of clarity  the order will be distinguished,  similarly to  Russo et al.  \cite{Russo2020},  using a combination of bars and tildes below the letter.
The equivalence between compact and indicial notation is then
\begin{equation}
\begin{gathered}
\tensf{a} \rightarrow  a_i\\
\tenss{a} \rightarrow  a_{ij} \\
\tenst{a} \rightarrow  a_{ijh}\\
\tensff{a} \rightarrow  a_{ijhk }\\
\end{gathered}
\nonumber
\end{equation} 
The tensor product between two first and two second order tensors is defined as
\begin{equation}
\begin{gathered}
\tenss{c}= \tensf{a} \otimes \tensf{b}  \rightarrow  c_{ij}=  a_i b_j \\
\tensff{c}=\tenss{a} \otimes \tenss{b} \rightarrow   c_{ijhk}  = a_{ij} b_{hk}
\end{gathered}
\nonumber
\end{equation}
The single contraction    (or scalar product) between a second order and first order tensors and between two second order tensors is 
\begin{equation}
\begin{gathered}
\tensf{c}=\tenss{a} \cdot \tensf{b} \rightarrow   c_i = a_{ij} b_j \\
\tenss{c}=\tenss{a} \cdot \tenss{b} \rightarrow  c_{ij}= a_{ik} b_{kj}
\end{gathered}
\nonumber
\end{equation}
whilst two \textbf{double} contractions  (or double scalar products) will be used
\begin{equation}
\begin{gathered}
c=\tenss{a} \colon \!  \tenss{b}= a_{ij} b_{ij} \\
c=\tenss{a} \cdotsd \tenss{b}= a_{ij} b_{ji}
\end{gathered}
\nonumber
\end{equation}
The gradient of  a tensor is 
\begin{equation}
 \tenss{c}= \tensf{a} \otimes \tensf{\nabla} \rightarrow c_{ij}= \der {a_i}{x_j} = a_{i,j}
 \nonumber
\end{equation}
The Levi-Civita \textit{pseudo}  or permutation tensor is indicated as 
\begin{equation}
\tenst{e} \rightarrow  e_{ijk}=
\begin{cases}
+1 \quad if \quad (i, j, k) \in \left\lbrace (1,2,3), (2,3,1), (3,1,2) \right\rbrace  \\
-1 \quad if \quad (i, j, k) \in \left\lbrace (1,3,2), (3,2,1), (2,1,3) \right\rbrace  \\
0  \quad \;\;   if \quad  i=j \quad  or \quad j=k  \quad or \quad i=k
\end{cases}
\nonumber
\end{equation}
and the following second and forth order unit tensors will be used 
\begin{equation}
\begin{gathered}
\tenss{I} \rightarrow \delta_{ij} \\
\tensff{I} \rightarrow  I_{ijhk}=\delta_{ih} \delta_{jk} \\
\bar{\tensff{I}} \rightarrow \bar{I}_{ijhk}=\delta_{ik} \delta_{jh} \\
\bar{\bar{\tensff{I}}} \rightarrow \bar{\bar{I}}_{ijhk}=\delta_{ij} \delta_{hk} \\
\bar{\tensff{I}}^{sym}= \frac{1}{2} \left( \tensff{I}+ \bar{\tensff{I}} \right) \rightarrow \bar{I}_{ijhk}^{sym} = \frac{1}{2} \left(  \delta_{ih} \delta_{jk} + \delta_{ik} \delta_{jh}   \right)  \\
\bar{\tensff{I}}^{skw}= \frac{1}{2} \left( \tensff{I}- \bar{\tensff{I}} \right) \rightarrow \bar{I}_{ijhk}^{skw} = \frac{1}{2} \left(  \delta_{ih} \delta_{jk} - \delta_{ik} \delta_{jh}   \right) 
\end{gathered}
\end{equation}
where the first and second are the so called \textit{identity} tensors such that
\begin{equation}
\begin{gathered}
\tenss{I} \cdot \tensf{a}= \tensf{a} \\
\tensff{I} \colon \tenss{a} = \tenss{a}
\end{gathered}
\end{equation}
whilst the remaining are such that 
\begin{equation}
\begin{gathered}
\bar{\tensff{I}} \colon \tenss{a}= \tenss{a}^T \\
\bar{\bar{\tensff{I}}} \colon \tenss{a}= \text{tr}(\tenss{a}) \tenss{I} \\
\tensff{I}^{sym} \colon \tenss{a}= \text{sym} \,  \tenss{a}\\
\tensff{I}^{skw} \colon \tenss{a}= \text{skew} \, \tenss{a}\\
\end{gathered}
\end{equation}
Finally the  unit deviatoric forth order tensor is defined as
\begin{equation}
\tensff{ \mathcal{I}}^d \rightarrow \mathcal{I}^d_{ijhk}=\delta_{ih} \delta_{jk} - \frac{1}{3} \delta_{hk} \delta_{ij}
\end{equation}
and is such that 
\begin{equation}
\tensff{\mathcal{I}}^d \colon \tenss{\varepsilon}= \tenss{e} \rightarrow \mathcal{I}^d_{ijhk} \varepsilon_{hk}=\left( \delta_{ih} \delta_{jk} - \frac{1}{3} \delta_{hk} \delta_{ij} \right) \varepsilon_{hk}= \varepsilon_{ij} - \frac{1}{3} \varepsilon_{rr} \delta_{ij}= e_{ij} 
\end{equation}

For the derivations presented in what follows we  will use the direct and inverse relationships between the vector which defines a rigid body rotation
\begin{equation}
\tensf{\omega}= \omega_i \hat{\tensf{e}}_i
\nonumber
\end{equation}
and its associated unit displacement gradient tensor $\tenss{\omega}$ 
\begin{equation}
\tenss{\omega}= -\tenst{e} \cdot \tensf{\omega} \qquad \rightarrow \qquad \omega_{ij}=-e_{ijk} \omega_k
\label{Eq_Tensor_Vector_Levi}
\end{equation}
and
\begin{equation}
\tensf{\omega}=-\frac{1}{2} \tenst{e} : \tenss{\omega}  \qquad \rightarrow \qquad \omega_k=-\frac{1}{2} e_{ijk} \omega_{ij}
\label{Eq_Vector_Tensor_Levi}
\end{equation}
where $\hat{\tensf{e}}_i$ with $i=1,2,3$ is an orthonormal basis.

\section{Cosserat continuum}
The so-called linear formulation of the general Cosserat  continuum (Cosserat \& Cosserat  \cite{Cosserat1909})  is adopted in this study.  In such a simplified version of the higher-order continuum developed by the two Cosserat brothers  the \textit{micro-volume} embedded in each point of the underlying material matrix  is assumed to be rigid (Mindlin  \cite{Mindlin1964}) and solidly connected to that point so that it moves together with it,  but it can also rigidly rotate relatively to the  material matrix.
Each   point  identified by the position vector 
\begin{equation}
\tensf{x}=x_i \bvh{e}_i
\nonumber
\end{equation}
 is thus characterized by six degrees of freedom (DoF) in the $3D$ space
\begin{equation}
\left\lbrace u_i,  \theta_i \right\rbrace
\nonumber
\end{equation}
where  $u_i$ and  $\theta_i$ are the components of the displacement and rotational  vector fields defined over the continuum
\begin{equation}
\begin{split}
\bv{u}=u_i \bvh{e}_i  \\
\bv{\theta}= \theta_i \bvh{e}_i
\end{split}
\nonumber
\end{equation}
Materials with a particulate nature,  such as granular media and soils,  when deformed not only experience displacements but also rotation of the individual particles.  The linear Cosserat continuum is hence inherently better suited than the classical Cauchy medium to reproduce their behaviour.

\subsection{Kinematics}
The unit variation of the  displacement vector in the neighbourhood of a material point is
\begin{equation}
\frac{d\tensf{u}}{dS}= \tensf{u} \otimes \tensf{\nabla} \cdot \frac{d\tensf{x}}{dS}=\tensf{u} \otimes \tensf{\nabla} \cdot \tensf{n} \quad \text{or} \quad \frac{du_i}{dS}= u_{i,j} \frac{dx_j}{dS}=  u_{i,j}  n_j
\nonumber
\end{equation}
where  $\tensf{n}$ is the unit directional vector pointing from  $\tensf{x}$  to  $\tensf{x}+d\tensf{x}$ and $dS$ is the distance between the two points
\begin{equation}
dS=\sqrt{ d\tensf{x} \cdot d\tensf{x} } \quad \text{or} \quad dS=\sqrt{dx_i dx_i}
\nonumber
\end{equation}
The gradient of the displacement vector field
\begin{equation}
\tensf{u} \otimes  \tensf{\nabla} \qquad	\text{or} \qquad u_{i,j}=\der{u_i}{x_j}
\nonumber
\end{equation}
can be split into its symmetric 
\begin{equation}
\tenss{\varepsilon}=\frac{1}{2} \left( \tensf{u} \otimes  \tensf{\nabla} +  \tensf{\nabla} \otimes \tensf{u}  \right) \qquad \text{or} \qquad \varepsilon_{ij}=\frac{1}{2}(u_{i,j}+u_{j,i}) 
\label{Eq_StrainDecompSymm}
\nonumber
\end{equation}
and skew-symmetric parts
\begin{equation}
\tenss{\omega}^m=\frac{1}{2} \left( \tensf{u} \otimes  \tensf{\nabla} -  \tensf{\nabla} \otimes \tensf{u}  \right) \qquad \text{or} \qquad  \omega^m_{ij}=\frac{1}{2}(u_{i,j}-u_{j,i})
\label{Eq_StrainDecompSkew}
\nonumber
\end{equation}
which represent the strain and the rigid body rotation tensors of the underlying material matrix, respectively. 
It should be noted that in the classical Cauchy continuum the  rotation  associated to the latter  tensor 
\begin{equation}
\tensf{\omega}^m=-\frac{1}{2} \tenst{e} : \tenss{\omega}^m   \qquad \textit{or}  \qquad \omega^m_k=-\frac{1}{2} e_{ijk} \omega^m_{ij}
\nonumber
\end{equation}
occurs freely  at constant energy,  hence   it is not associated to a stress variation. 
However, in the case of the Cosserat medium such a rotation does not necessarily occur freely, without providing any work,  as it  interacts with the rotation $\bv{\theta}$   prescribed to  the \textit{micro-volumes} which are themselves embedded in the matrix  at each point  $\bv{x}$. 
If the  rotation  of the embedded \textit{micro-volumes}  is left free to follow that of  the underlying matrix,  then no work needs to be provided and no stress variation occurs.  However if the matrix rotation $\bv{\omega}^m$ is larger than the  rotations $\bv{\theta}$ prescribed to the \textit{micro-volumes},  then work needs to be provided to rotate the matrix and  additional stress will develop.

This implies that the   Cosserat strain tensor must be   defined as
\begin{equation}
\tenss{\gamma}= \tensf{u} \otimes \tensf{\nabla} - \tenss{\theta} \qquad \text{or} \qquad \gamma_{ij}=u_{i,j}-\theta_{ij}
\nonumber
\end{equation}
or equivalently
\begin{equation}
\tenss{\gamma}= \tenss{\varepsilon} + \tenss{\omega}  \qquad \textit{or} \qquad \gamma_{ij}=\varepsilon_{ij} + \omega_{ij}
\label{gamma_def}
\end{equation}
where $\tenss{\omega}$ is usually referred to  as the \textit{relative rotation} and defined as
\begin{equation}
 \tenss{\omega}  =  \tenss{\omega}^m-\tenss{\theta}  \qquad \textit{or} \qquad  \omega_{ij}= \omega^m_{ij} - \theta_{ij}
 \nonumber
\end{equation}
or
\begin{equation}
 \tenss{\omega}  = - \tenst{e}\cdot(\tensf{\omega}^m-\tensf{\theta})   \qquad \textit{or} \qquad  \omega_{ij}= - e_{ijk} (\omega^m_k - \theta_k)
 \nonumber
\end{equation}
and represents the relative rotation of the material matrix with respect to that  prescribed to the  \textit{micro-volumes} embedded at each location  $\tensf{x}$.

\begin{figure}[tb]
\centering
\includegraphics[width=\textwidth]{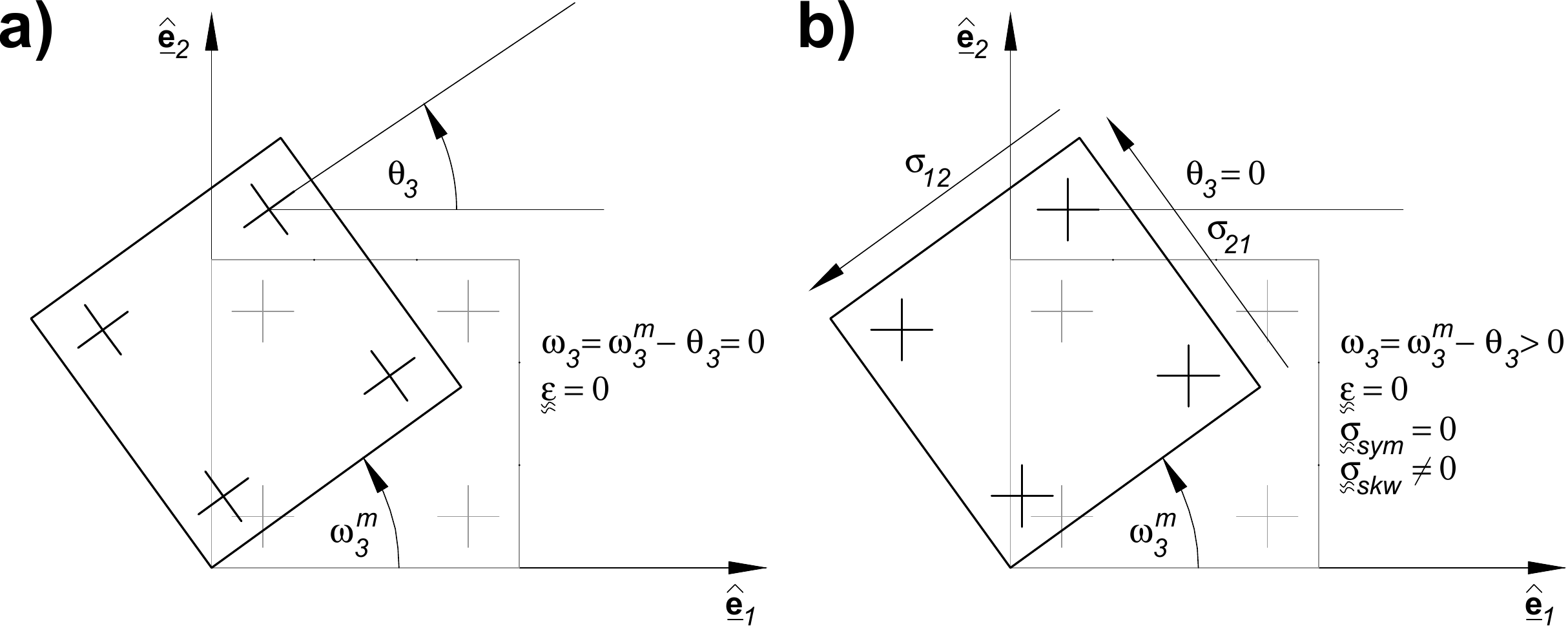}
\caption{Role played by the relative rotation tensor into the Cosserat strain and stress definitions in a plane  strain exemplification of the Cosserat continuum.}
\label{Fig_Kinematics}
\end{figure}

Fig. \ref{Fig_Kinematics} shows a plane strain exemplification of the role played by the relative rotation into the Cosserat strain definition.  The number of DoF reduces to three 
\begin{equation}
\begin{gathered}
\tensf{u} =u_1 \hat{\tensf{e}}_1 +u_2 \hat{\tensf{e}}_2 \\
\tensf{\theta} =\theta_3 \hat{\tensf{e}}_3
\end{gathered}
\nonumber
\end{equation}
In Fig. \ref{Fig_Kinematics}.a  the  strain tensor of the material matrix  $\tenss{\varepsilon}$ is nil, whilst its rotation and the rotation of the embedded \textit{micro-volumes} have the same value, so that the relative rotation  $\tenss{\omega}$ vanishes.  No stress state acts on such a Cosserat infinitesimal element. 
In Fig. \ref{Fig_Kinematics}.b  the strain tensor of the  material  matrix is still nil and its rotation is $\tenss{\omega}^m>0$, whilst the rotation of the embedded \textit{micro-volume} is set to zero, so that a positive relative rotation $\tenss{\omega}$ applies.  Whilst the stress state associated to the deformation $\tenss{\varepsilon}$ of the material matrix is nil, the anti-clockwise and positive rotation of the matrix $\tensf{\omega}^m$  requires an anti-clockwise tangential stress state on the material element.
It should be noted that whilst the stress tensor associated to the symmetric part of the strain tensor $\tenss{\gamma}$ is also symmetric, that of the skew-symmetric part of $\tenss{\gamma}$ is also skew-symmetric.

In addition to the strain tensor a curvature or wryness tensor characterizes the Cosserat continuum
\begin{equation}
\tenss{\chi}= \tensf{\theta} \otimes \tensf{\nabla} \qquad \text{or} \qquad \chi_{ij}=\theta_{i,j}
\label{chi_def}
\end{equation}
which accounts for the different prescribed rotations of infinitesimally close  \textit{micro-volumes}.  
It should be noted that  the principal diagonal  $\chi_{i,i}$ and the off-diagonal  $\chi_{i,j}$ terms   of this tensor,  i.e.  its  spherical  and deviatoric components,  represent torsional and bending  curvatures, respectively.

\subsection{Statics}
Let $\tenss{\sigma}$ and $\tenss{\mu}$  indicate the stress and couple-stress tensors conjugated in the meaning of Hill \cite{Hill1968b} to the strain $\tenss{\gamma}$  and wryness $\tenss{\chi}$ tensors  introduced in Eqs.  \eqref{gamma_def} and \eqref{chi_def}.
The rate of the internal  work of deformation  over the volume $\Omega$ is 
\begin{equation}
\dot{W}_{int}=\int_{\Omega} \left( \tenss{\sigma} : \dot{\tenss{\gamma}} + \tenss{\mu}:\dot{\tenss{\chi}} \right) \, d\Omega
\label{InternalWork_gen}
\end{equation}
However this can  be reformulated  by exploiting the decomposition of a general second order tensor into its symmetric and skew-symmetric parts, which for the Cosserat stress and strain tensors gives
\begin{equation}
\begin{gathered}
\tenss{\sigma}=\text{sym} \,  \tenss{\sigma}+\text{skew} \, \tenss{\sigma} = \tenss{\sigma}_{sym}+\tenss{\sigma}_{skw} \\
\tenss{\gamma}=\text{sym} \, \tenss{\gamma} +\text{skew} \, \tenss{\gamma}= \tenss{\varepsilon} + \tenss{\omega}
\end{gathered}
\end{equation}
Since the double contraction of a symmetric and a skew-symmetric tensor is nil,    the term $\tenss{\sigma}:\dot{\tenss{\gamma}}$ can be rewritten as
\begin{equation}
\tenss{\sigma}:\dot{\tenss{\gamma}}=\left( \text{sym}\tenss{\sigma}+\text{skew}\tenss{\sigma}  \right) \colon \left( \text{sym}\dot{\tenss{\gamma}} +\text{skew}\dot{\tenss{\gamma}} \right)= \tenss{\sigma}_{sym} \colon \dot{\tenss{\varepsilon}} + \tenss{\sigma}_{skw} \colon \dot{\tenss{\omega}}
\end{equation}
and Eq. \eqref{InternalWork_gen} then becomes
\begin{equation}
\dot{W}_{int}=\int_{\Omega} \left( \tenss{\sigma}_{sym} : \dot{\tenss{\varepsilon}} + \tenss{\sigma}_{skw} : \dot{\tenss{\omega}} + \tenss{\mu}:\dot{\tenss{\chi}} \right) \, d\Omega
\label{Eq_InternalWork_symskew}
\end{equation}
which indicates that the symmetric part of the stress tensor performs work only for the classical Cauchy strain tensor whilst the skew-symmetric part of the same tensor performs work for the \textit{relative-rotation} tensor only.

By applying the  Green's theorem to  $\tenss{\sigma}_{sym} \colon \dot{\tenss{\varepsilon}}$ ,   $\tenss{\sigma}_{skw} \colon \dot{\tenss{\omega}}$  and to $\tenss{\mu} \colon \dot{\tenss{\chi}}$ one obtains
\begin{equation}
\begin{split}
\int_\Omega \tenss{\sigma}_{sym} \colon \dot{\tenss{\varepsilon}} \, d\Omega= \int_{\partial\Omega} \dot{\tensf{u}} \cdot \tenss{\sigma}_{sym} \cdot \tensf{n} \, dS- \int_\Omega \dot{\tensf{u}} \cdot \tenss{\sigma}_{sym}\cdot \tensf{\nabla} \, d\Omega \\
\int_\Omega \tenss{\sigma}_{skw} \colon \dot{\tenss{\omega}} \, d\Omega= 
 \int_{\partial\Omega} \dot{\tensf{u}} \cdot \tenss{\sigma}_{skw} \cdot \tensf{n} \, dS- \int_\Omega \dot{\tensf{u}} \cdot \tenss{\sigma}_{skw}\cdot \tensf{\nabla} \, d\Omega - \int_\Omega 2 \tensf{\sigma}_{skw} \cdot \dot{\tensf{\theta}} \,  d\Omega \\
 \int_\Omega \tenss{\mu} \colon \dot{\tenss{\chi}} \,  d\Omega= 
 \int_{\partial\Omega} \dot{\tensf{\theta}} \cdot \tenss{\mu} \cdot \tensf{n} \, dS - \int_\Omega  \dot{\tensf{\theta}} \cdot \tenss{\mu}  \cdot \tensf{\nabla} \, d\Omega
 \end{split}
\nonumber
\end{equation}
where $\bv{n}$ is the outer normal to the continuum,  $\partial\Omega$ is the surface enclosing the continuum and Eq. \eqref{Eq_Tensor_Vector_Levi} has been used to transform the skew stress and \textit{micro-volumes} rotation tensors into their  vector counterparts.

Finally the substitution of these results into Eq. \eqref{Eq_InternalWork_symskew} yields 
\begin{equation}
\begin{aligned}
\dot{W}_{int}=
  \int_{\partial\Omega} \dot{\bv{u}} \cdot \tenss{\sigma} \cdot  \bv{n} \, dS -  \int_\Omega\dot{\bv{u}} \cdot \tenss{\sigma} \cdot \bv{\nabla} \, d\Omega - \int_\Omega 2 \bv{\sigma}_{skw} \cdot  \dot{\bv{\theta}} \, d\Omega + \\
\int_{\partial\Omega} \dot{\bv{\theta}} \cdot \tenss{\mu} \cdot \bv{n} \,dS - \int_\Omega \dot{\bv{\theta}} \cdot \tenss{\mu} \cdot \bv{\nabla} \,  d\Omega\\
\end{aligned}
\end{equation}

The rate of work of deformation of the external forces over the volume $\Omega$ and bounding surface $\partial{ \Omega}$ of the continuum  is
\begin{equation}
\dot{W}_{ext}=\int_{\Omega} \left( \tensf{f}  \cdot \dot{\tensf{u}} + \tensf{c} \cdot \dot{\tensf{\theta}} \right) \, d\Omega + \int_{\partial\Omega} \left( \hat{\tensf{t}} \cdot \dot{\tensf{u}} +\hat{ \tensf{m}} \cdot \dot{\tensf{\theta}}  \right) \,dS
\label{ExternalWork_gen}
\end{equation}
where $\tensf{f}$ and $\tensf{c}$ are the body forces and couples per unit volume and  $\hat{\tensf{t}}$ and $\hat{\tensf{m}}$ are the surface tractions and couples per unit surface.

The balance equation of the internal and external work of deformation  in static applications  can then be written as 

\begin{align}
\begin{split}
\dot{W}_{ext} - \dot{W}_{int} =
& \int_{\partial\Omega} \dot{\tensf{u}} \cdot \left( \hat{\tensf{t}} - \tenss{\sigma} \cdot \tensf{n} \right) \, dS + \int_{\partial\Omega} \dot{\tensf{\theta}} \cdot \left( \hat{\tensf{m}} - \tenss{\mu} \cdot \tensf{n} \right) \, dS + \\
&  \int_\Omega \dot{\tensf{u}} \cdot \left( \tenss{\sigma} \cdot \tensf{\nabla} + \tensf{f} \right) \,  d\Omega + 
  \int_\Omega \dot{\tensf{\theta}} \cdot \left( 2 \tensf{\sigma}_{skw} + \tensf{c} + \tenss{\mu} \cdot \bv{\nabla} \right) = 0
\end{split}
\end{align}
from which the strong form of the balance equations is retrieved
\begin{equation}
\begin{gathered}
 \tenss{\sigma} \cdot \bv{\nabla} + \bv{f}=0 \\
 \tenss{\mu} \cdot \bv{\nabla} + 2 \bv{\sigma}_{skw} + \bv{c} =0
\end{gathered}
\label{Eq_Balance_Strong}
\end{equation}
and on the surface of the continuum
\begin{equation}
\begin{gathered}
\tenss{\sigma} \cdot \tensf{n}=\hat{\tensf{t}} \\
 \tenss{\mu} \cdot \tensf{n}=   \hat{\tensf{m}}
\end{gathered}
\end{equation}
The second of Eqs \eqref{Eq_Balance_Strong} shows that due to  the couple-stress tensor the stress tensor is not symmetric and that its skew component is 
\begin{equation}
\tenss{\sigma}_{skw}=  \tenst{e} \cdot \frac{1}{2} \left(  \tenss{\mu} \cdot \tensf{\nabla} + \tensf{c} \right)
\end{equation}

\begin{figure}[tb]
\centering
\includegraphics[width=0.7\textwidth]{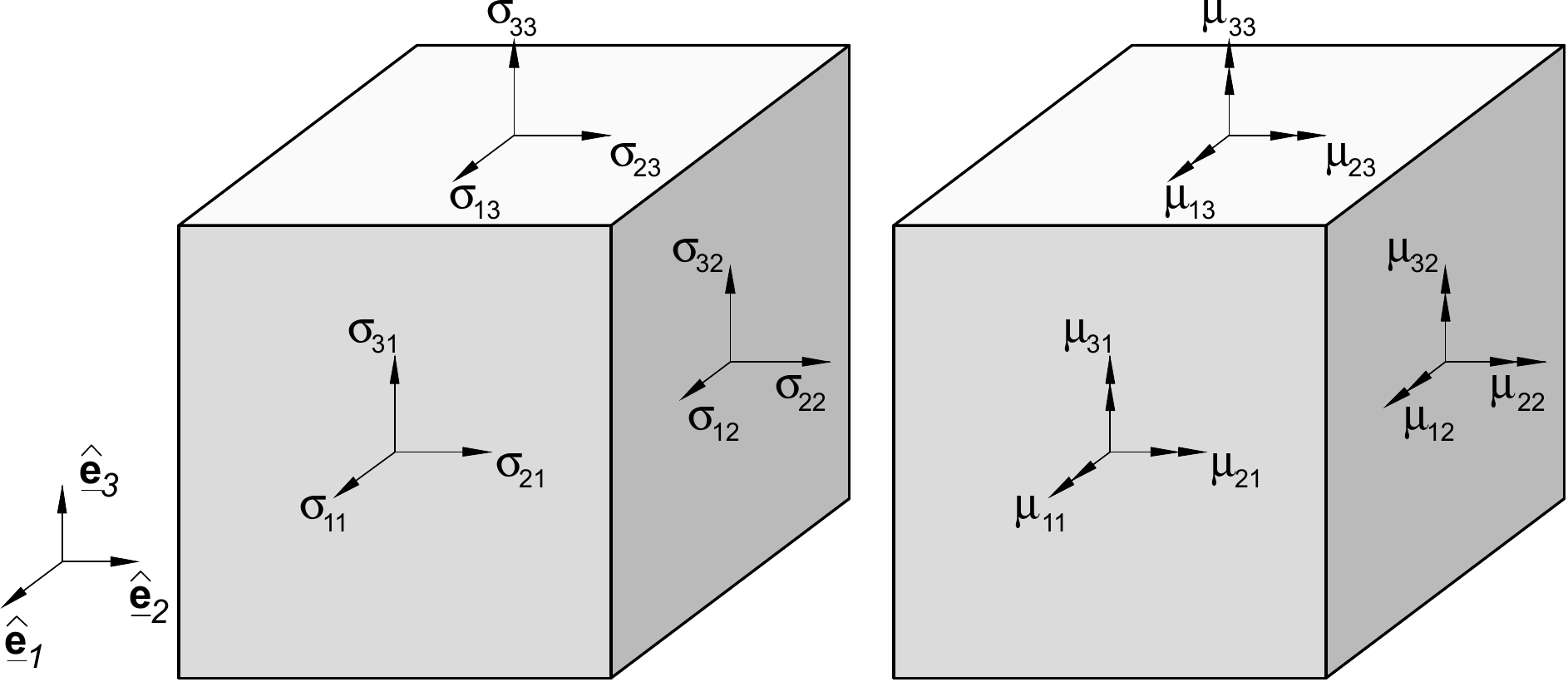}
\caption{Stress and couple-stress tensors for the Cosserat continuum.}
\label{Fig_SignConvention}
\end{figure}

The sign convention for both stress and couple-stress tensors is shown in Fig \ref{Fig_SignConvention}. It should be noted that the second subscript indicates the normal to the plane.

\section{Elastic behaviour}
It is assumed that the reversible behaviour of the Cosserat continuum is described by means of a linear hyper-elastic model.
The work of deformation  is then stored in each material point  in the form of an  elastic potential energy  $\Psi$ defined as
\begin{equation}
\Psi = \int d\Psi=  \int \left(\tenss{\sigma} \colon \dot{\tenss{\gamma}}^e + \tenss{\mu} \colon \dot{\tenss{\chi}}^e  \right)
\label{Eq_ElasticStrainEnergyPotential_Gen}
\end{equation}
where $\tenss{\gamma}^e$ and $\tenss{\chi}^e$ indicate the elastic components of the respective tensors.
However,  the decomposition of all stress and strain tensors into their symmetric and skew-symmetric parts allows a simplification based on the observation that  the double contraction of a symmetric and skew-symmetric tensor is ni. 
Eq. \eqref{Eq_ElasticStrainEnergyPotential_Gen} can  hence be rewritten as
\begin{equation}
\Psi= \int d\Psi=  \int \left(\tenss{\sigma}_{sym} \colon \dot{\tenss{\varepsilon}}^e + \tenss{\sigma}_{skw} \colon \dot{\tenss{\omega}}^e + \tenss{\mu}_{sym} \colon \dot{\tenss{\chi}}_{sym}^e  +  \tenss{\mu}_{skw} \colon \dot{\tenss{\chi}}_{skw}	^e  \right)
\label{Eq_ElasticStrainEnergyPotential_GenSymSkew}
\end{equation}
and since   $d\Psi$ is a perfect differential 
\begin{equation}
\begin{gathered}
\tenss{\sigma}_{sym}= \frac{\partial \Psi}{ \partial \tenss{\varepsilon}^e}  \qquad  \qquad  
\tenss{\sigma}_{skw}= \frac{\partial \Psi}{ \partial \tenss{\omega}^e} \\
\tenss{\mu}_{sym}= \frac{\partial \Psi}{ \partial \tenss{\chi}_{sym}^e} \qquad  \qquad
\tenss{\mu}_{skw}= \frac{\partial \Psi}{ \partial \tenss{\chi}_{skw}^e}
\end{gathered}
\end{equation}

In the case of linear elastic behaviour the energy potential is a quadratic function given by 
\begin{multline}
\Psi\left( \tenss{\varepsilon}^e ,  \tenss{\omega}^e, \tenss{\chi}_{sym}^e ,  \tenss{\chi}_{skw}^e \right)= \\
   \frac{1}{2} \tenss{\varepsilon}^e \colon \frac{\partial^2 \Psi}{\partial \tenss{\varepsilon}^e \otimes  \partial \tenss{\varepsilon}^e} \colon \tenss{\varepsilon}^e + 
\frac{1}{2} \tenss{\omega}^e \colon \frac{\partial^2 \Psi}{\partial \tenss{\omega}^e \otimes \partial{\tenss{\omega}^e}} \colon \tenss{\omega}^e + \\
 \frac{1}{2} \tenss{\chi}_{sym}^e \colon \frac{\partial ^2 \Psi}{\partial \tenss{\chi}_{sym}^e \otimes \partial \tenss{{\chi}}_{sym}^e} \colon \tenss{\chi}_{sym}^e + \frac{1}{2} \tenss{\chi}_{skw}^e \colon \frac{\partial ^2 \Psi}{\partial \tenss{\chi}_{skw}^e \otimes \partial \tenss{{\chi}}_{skw}^e} \colon \tenss{\chi}_{skw}^e 
\end{multline}
and  the constant forth order isotropic symmetric and skew symmetric stiffness tensors are
\begin{equation}
\begin{gathered}
\tensff{D}_\varepsilon=\frac{\partial^2 \Psi}{\partial \tenss{\varepsilon}^e \otimes  \partial \tenss{\varepsilon}^e} = \left( K -\frac{2}{3} G \right) \bar{\bar{\tensff{I}}} + 2 G \tensff{I}^{sym}  \\
\tensff{D}_\omega= \frac{\partial^2 \Psi}{\partial \tenss{\omega}^e \otimes \partial{\tenss{\omega}^e}} = 2 \, G_c \tensff{I}^{skw} \\
\tensff{D}_{\chi_{sym}}=\frac{\partial^2 \Psi}{\partial \tenss{\chi}_{sym}^e \otimes \partial {\tenss{\chi}_{sym}}} = T \bar{\bar{\tensff{I}}} + 2 \, B \tensff{I}^{sym}\\
\tensff{D}_{\chi_{skw}}=\frac{\partial^2 \Psi}{\partial \tenss{\chi}_{skw}^e \otimes \partial {\tenss{\chi}_{skw}}}= 2 \, B_c \tensff{I}^{skw}
\end{gathered}
\label{Eq_ElasticStiffnessTensors}
\end{equation}
In these equations $K$ is the bulk modulus,  hence $K-\frac{2}{3}G$  and $G$ are the usual Lame's constants ( clearly associated to the symmetric of the tensors only) whilst  $G_c$ is an additional constitutive parameter also with dimension of a force per unit area.   $T$, $B$ and $B_c$ are additional parameters which account for the torsional and bending stiffness with dimension of a moment per unit length.  

The linear elastic stress-strain relationships can then  be formulated as
\begin{equation}
\begin{gathered}
\tenss{\sigma}= \left(K-\frac{2}{3}G\right) \text{tr}{\tenss{\varepsilon}^e} \tenss{I} + 2 G \tenss{\varepsilon}^e + 2 G_c \tenss{\omega}^e \\
\tenss{\mu}= T \text{tr}\tenss{\chi}^e \tenss{I} + 2 B \tenss{\chi}^e_{sym} + 2 B_c \tenss{\chi}^e_{skw}
\end{gathered}
\label{Eq_ElasticLawNoDecomp}
\end{equation}

The latter can be further rewritten considering the spherical and deviatoric components of the stress  and couple-stress tensors and of the strain and wryness tensors
\begin{equation}
\begin{gathered}
 \tenss{s}_{sym}=\tenss{\sigma}_{sym} - \frac{1}{3}\text{tr} \tenss{\sigma} \tenss{I}  \\
 \tenss{m}_{sym}=\tenss{\mu}_{sym} - \frac{1}{3} \text{tr} \tenss{\mu} \tenss{I}   \\
 \tenss{e}_{sym}=\tenss{\varepsilon}-\frac{1}{3} \text{tr} \tenss{\varepsilon} \tenss{I} \\
 \tenss{g}_{sym}=\tenss{\chi}-\frac{1}{3} \text{tr} \tenss{\chi} \tenss{I} \\
\end{gathered}
\label{Eq_DeviatorcSphericalTensors}
\end{equation}
whilst the skew-symmetric part of all four tensors is clearly deviatoric only.
The elastic constitutive equations can then be reformulated as
\begin{equation}
\begin{gathered}
\tenss{\sigma}= K  \text{tr}{\tenss{\varepsilon}^e} \tenss{I} + 2 G \tenss{e}^e + 2 G_c \tenss{\omega}^e \\
\tenss{\mu}= K_c \text{tr}\tenss{\chi}^e \tenss{I} + 2 B \tenss{g}^e_{sym} + 2 B_c \tenss{g}^e_{skw}
\end{gathered}
\label{Eq_StressStrainSpherDev}
\end{equation}
where to simplify the notation
\begin{equation}
K_c=(T + \frac{2}{3} B)
\end{equation}

\section{Equivalent von Mises stress}

For the formulation of the yield function a definition of the \textit{ equivalent von Mises stress} is required for the Cosserat continuum.
The  von Mises  \cite{vonMises1913} criterion assumes that yield in a Cauchy continuum occurs when the second invariant $J_2=\frac{1}{2} \tenss{s}_{sym} \colon \tenss{s}_{sym}$ of the deviatoric stress tensor  attains a critical value.
According to the  physical interpretation of  that criterion given by Hencky \cite{Hencky1924},  yield in a specific material occurs when the  energy of the  linear elastic distortion reaches a critical value.  This can be obtained by subtracting from Eq. \eqref{Eq_ElasticStrainEnergyPotential_GenSymSkew} the energy associated to the volumetric deformation
\begin{equation}
\Psi^e_D=\frac{1}{2}\left[ \tenss{s}_{sym} \colon \tenss{\varepsilon}^e + \tenss{s}_{skw} \colon \tenss{\omega}^e + \tenss{m}_{sym} \colon \tenss{g}_{sym}^e + \tenss{m}_{skw} \colon \tenss{g}_{skw}^e + \frac{1}{3} \text{tr} \tenss{\mu}  \text{tr} \tenss{\chi}^e \right]
\end{equation}
where the decomposition given in Eqs \eqref{Eq_DeviatorcSphericalTensors}  of the stress, strain, couple stress and curvature  tensors into their respective deviatoric and spherical parts has been used. 
It should be noted that whilst the spherical component of the strain  tensor is indeed not associated to distortion,  the principal diagonal terms of the curvature tensors account for the rate of variation along one axis of the  component along the same axis of  the rotation vector of the  \textit{micro-volumes}, i.e. they define a torsional curvature (e.g.  Russo et al. \cite{Russo2020}).  Since those terms are clearly related  to a distortion,  the energy associated to the spherical component of the curvature and couple-stress tensors needs to be included   in the definition of the energy of distortion.  

After the substitution of   the  linear elastic constitutive relationship of Eqs.  \eqref{Eq_StressStrainSpherDev} the energy becomes
\begin{equation}
\begin{split}
\Psi^e_D=\frac{1}{4 G} \left[ \tenss{s}_{sym} \colon \tenss{s}_{sym} + \frac{G}{G_c} \tenss{s}_{skw} \colon \tenss{s}_{skw} + \frac{G}{B} \tenss{m}_{sym} \colon \tenss{m}_{sym} + \right. \\
\left.  \frac{G}{B_c} \tenss{m}_{skw} \colon \tenss{m}_{skw} + \frac{2 G}{ K_c} \frac{ \text{tr}^2\tenss{\mu}}{9} \right]
\end{split}
\end{equation}
If the energy associated to the linear elastic distortion is the quantity responsible for yielding,  then by imposing the equivalence with that of a bar of a Cauchy material subjected to a uni-axial stress state $q$  
\begin{equation}
\Psi^e_D=\frac{1}{2} \frac{q^2}{3 G}
\end{equation}
the \textit{equivalent von Mises stress} is obtained
\begin{equation}
\begin{split}
q= \left\lbrace \frac{3}{2} \left[ \tenss{s}_{sym} \colon \tenss{s}_{sym} + \frac{G}{G_c} \tenss{s}_{skw} \colon \tenss{s}_{skw} + \frac{G}{B} \tenss{m}_{sym} \colon \tenss{m}_{sym} +\frac{G}{B_c} \tenss{m}_{skw} \colon \tenss{m}_{skw} \right] \right.  \\
\left.   + \frac{2 G}{K_c} \frac{ \text{tr}^2\tenss{\mu}}{9}  \right\rbrace ^{\frac{1}{2}}
\end{split}
\label{Eq_vonMisesEqStressCosserat}
\end{equation}
For dimensional correctness it is clear that three quantities
\begin{equation}
\begin{gathered}
\ell_1=\sqrt{ \frac{B}{G}} \\
\ell_2=\sqrt{\frac{B_c}{G}}\\
\ell_3=\sqrt{\frac{K_c}{2G}}
\end{gathered}
\end{equation}
have the dimension of a length and can be considered as three	 \textit{characteristic lengths} of the Cosserat material.
It should be noted that in plane strain problems  the third characteristic length does not play any role, since the spherical component of the curvature tensor vanishes.

This result is coherent with the theory of representations which  states that an isotropic scalar valued function of a finite number of symmetric $\tenss{A}_i$ and skew-symmetric $\tenss{W}_i$ second order tensors has a  \textit{complete} and \textit{irreducible} representation $f$  in terms of a number of invariants  expressed as the trace of the scalar product  of those tensors.  A list of those invariants for different situations can be found in e.g.  Smith and  Smith, \cite{Smith1981},  Zheng \cite{Zheng1994}.
In the two simplest cases where the arguments of the function are either a single symmetric tensor or a symmetric and a skew-symmetric tensors, the \textit{complete} and \textit{irreducible} representations are 
\begin{equation}
f=f( \text{tr} \tenss{A}, \text{tr}\tenss{A}^2, \text{tr} \tenss{A}^3)
\end{equation}
and
\begin{equation}
f= f\left( \text{tr} \tenss{A},  \text{tr}\tenss{A}^2,  \text{tr} \tenss{A}^3, \text{tr} \tenss{W}^2,  \text{tr} \tenss{A} \cdot \tenss{W}^2,  \text{tr} \tenss{A}^2 \cdot \tenss{W}^2,  \text{tr} \tenss{A}^2 \cdot \tenss{W}^2 \cdot \tenss{A} \cdot \tenss{W}  \right)
\label{Eq_GeneralRepresentationSymSkew}
\end{equation}
respectively.
The former expression is often used in Cauchy continua to formulate  the isotropic scalar-valued functions which define yield and failure criteria  (e.g.  von Mises and Drucker-Prager)  and  strain energy potentials  (e.g.  Houlsby \cite{Houlsby1985},    Lagioia and Panteghini \cite{Lagioia2019}), often dropping  one or two  of the three invariants,  resulting in \textit{incomplete} formulations.  For example the third invariant is dropped in the   Drucker-Prager criterion,  thus loosing  the dependency  on the  relative magnitude of the intermediate principal stress $\sigma_{II}$ with respect to $\sigma_I$ and $\sigma_{III}$, where $\sigma_I\geq\sigma_{II}\ge \sigma_{III}$.
Similarly if in  Eq.  \eqref{Eq_GeneralRepresentationSymSkew} some of the invariants are dropped,  retaining only  $ \text{tr} \tenss{A}^2$,  $ \text{tr} \tenss{W}^2$,  then a general expression for the isotropic scalar-valued function of the equivalent von Mises stress of  Eq. \eqref{Eq_vonMisesEqStressCosserat} is obtained ( to be precise two skew-symmetric tensors are used in that case). 

Finally it should be noted that since 
\begin{equation}
\begin{gathered}
a_3 \text{tr} \tenss{\sigma}_{sym}^2 +a_4 \text{tr} \tenss{\sigma}^2_{skw}  \\
a_1 \tenss{\sigma} \colon \tenss{\sigma} + a_2 \tenss{\sigma} \colon \tenss{\sigma}^T
\end{gathered}
\end{equation}
are equivalent, provided that
\begin{equation}
\begin{gathered}
a_1=\frac{a_3-a_4}{2} \\
a_2=\frac{a_3+a_4}{2}
\end{gathered}
\end{equation}
 for  plane strain conditions,  Eq. \eqref{Eq_vonMisesEqStressCosserat} can also be formulated as
\begin{equation}
\begin{split}
q= \left\lbrace \frac{3}{2} \left[  \frac{G_c-G}{2G_c} \tenss{s} \colon \tenss{s}+ \frac{G_c+G}{2G_c} \tenss{s} \colon \tenss{s}^T + \frac{G}{B} \frac{B_c-B}{2B_c}\tenss{m} \colon \tenss{m} \right.  \right. \\
\left. \left. +\frac{G}{B}\frac{B_c+B}{2B_c} \tenss{m}\colon \tenss{m}^T \right]  \right\rbrace ^{\frac{1}{2}}
\end{split}
\label{Eq_vonMisesEqStressCosseratTranspose}
\end{equation}
which coincides with that  proposed  by Muehlhaus and Vardulakis \cite{Muehlhaus1987}.
It is interesting to observe that those Authors derived their definition of $q$ using only micromechanical and averaging considerations, and yet it is equivalent to that obtained using either the theory of representations or Hencky's energy approach.

\section{Constitutive model}
We will confine our work to isotropic materials.
Constitutive models,  even with  such a frequent and simplifying assumption,  are usually formulated and numerically integrated  in terms of components of the stress and strain tensors.  However it has been recently shown by Panteghini and Lagioia \cite{PL2014} \cite{PL2018} that if  stress and strain invariants are adopted,  an extremely efficient integration algorithm can be written which requires the solution of one  equation in one unknown only, rather than a system  of seven by seven equations and unknowns as in the standard integration. This not only results in considerable reduction of the machine run-time but also in a more robust scheme.

The approach described by Panteghini and Lagioia  is here further extended to  formulate and integrate an elasto-plastic constitutive model for a Cosserat continuum. 
The recoverable behaviour is assumed to be hyper-elastic and linear as described previously, whilst the irrecoverable one is defined within the  theory of plasticity with isotropic hardening/softening.  

The adopted framework is that typically  used for elastic-perfect plastic materials,  where the yield and the plastic potential surfaces are described by classical yield/failure criteria.  However an hardening rule is also introduced which governs the evolution of the intercept of the yield surface with the equivalent von Mises stress axis.

The  isotropic scalar-valued function which defines the yield and plastic potential surfaces needs to be defined entirely in terms of stress invariants
\begin{equation}
f=f(p, q, \theta_s)
\nonumber
\end{equation}
where $p$, $q$ and $\theta_s$ for a Cauchy continuum are the mean pressure,   equivalent von Mises stress (also known as  \textit{deviatoric stress}) and the Lode's angle of the stress tensor.

This expression can also  be employed  for a Cosserat continuum,  if   the definition of $p$ is kept unchanged,  i.e.    $p= \frac{\text{tr} \tenss{\sigma}}{3}$,  whilst the deviatoric stress is  substituted with that given by  Eq. \eqref{Eq_vonMisesEqStressCosserat} and the Lode's angle is that of the symmetric part only of the deviatoric stress $\tenss{s}_{sym}$ tensor
\begin{equation}
\theta_s= \frac{1}{3} \arcsin \left(  -\frac{27}{2} \frac{\text{det} \tenss{s}_{sym}}{q_s^3}\right)
\nonumber
\end{equation}
where 
\begin{equation}
q_s= \sqrt{\frac{3}{2} \tenss{s}_{sym} \colon \tenss{s}_{sym}}
\end{equation}

Following  Panteghini and Lagioia \cite{PL2018}  considerable advantages are achieved if a slightly different structure is adopted 
\begin{equation}
f=f\left( p, q \Gamma({\theta}_s) \right)
\label{Eq_YieldFunctionGeneral}
\end{equation}
in which  $\Gamma(\theta_s)$ is the reciprocal of the  function which defines the shape of the yield and plastic potential surfaces  in the deviatoric plane,  and the Lode's angle  $\theta_s$ is restricted in the interval  $\theta_s \in \left[ -\frac{\pi}{6}, \frac{\pi}{6}, \right]$. 
It should be noted that in triaxial compression conditions, associated to $\theta_s=-\frac{\pi}{6}$,   $\Gamma$ evaluates to unity, whilst  in triaxial extension conditions, characterized by $\theta_s=\frac{\pi}{6}$,  the value of $\Gamma$ is larger than unity  and results $\Gamma\left( \frac{\pi}{6}\right)= \frac{q_y(-\frac{\pi}{6})} {q_y ( \frac{\pi}{6})}$,  where the subscript $y$ indicate  yield in the case of metals and failure in the case of soils. 
The role of $\Gamma(\theta_s)$ is that of modulating the value of the equivalent von Mises stress $q$ so that the meridional sections of the yield and plastic potential surfaces are accordingly reduced or expanded.  

The yield and the plastic potential surfaces are defined by a function recently proposed by Lagioia and Panteghini \cite{LP2016}. That function was mathematically demonstrated to define classical yield and failure criteria for metals and soils,  i.e.  von Mises,  Tresca,  Drucker-Prager, Matsuoka-Nakai,  Lade-Duncan and Mohr-Coulomb.  For the last failure criterion parameters were also retrieved which result in smooth  versions circumscribed to and approximately inscribed in the original one.   The set of published parameters is likely to be not conclusive and other classical criteria can possibly be included.  As an example,   recently Lester and Sloan \cite{Lester2017} retrieved parameters for a rigorously inscribed rounded Mohr-Coulomb and named the Lagioia and Panteghini formulation  \textit{generalized classical}  yield function (GC). 

\begin{figure}[tb]
\centering
\includegraphics[width=\textwidth]{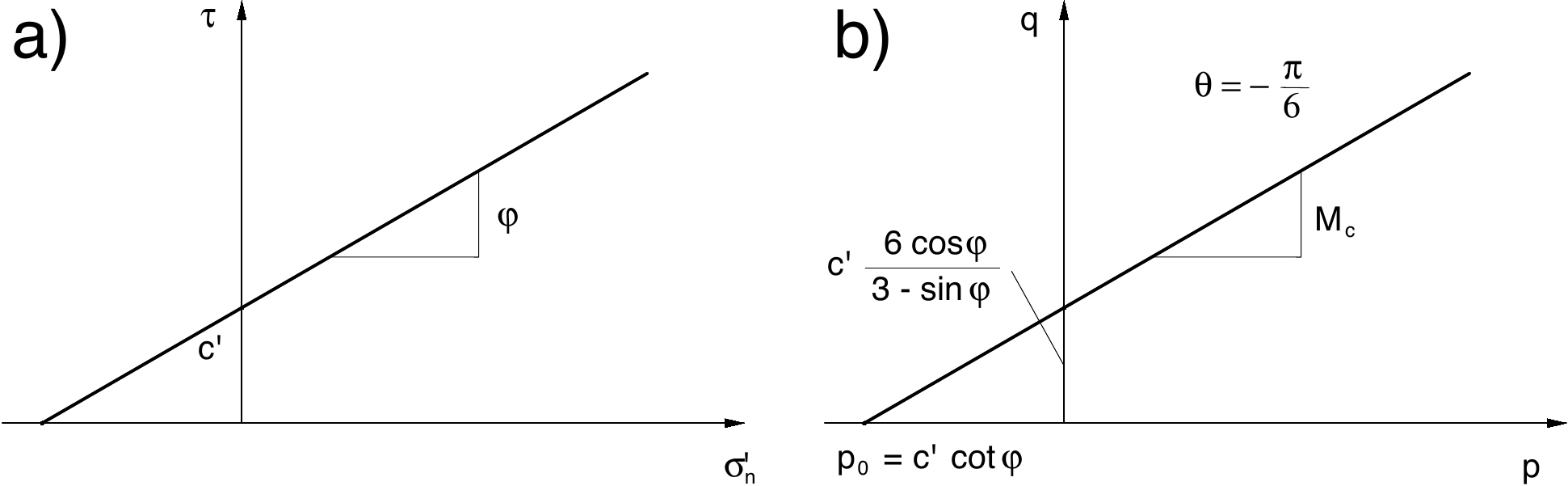}
\caption{Significance of the \textit{cohesion} term in the failure criterion.}
\label{Fig_CohesionTerm}
\end{figure}

Even if Eq. \ref{Eq_YieldFunctionGeneral} can be used to define yield surfaces with any meridional section,  for for  those criteria  it is linear and 
 Eq. \ref{Eq_YieldFunctionGeneral} particularizes to the Generalised Classical yield function, which with the  sign convention adopted in this paper  (tensile stresses are positive) becomes
\begin{equation}
f(p,q,\theta_s)= q \,  \Gamma \left(\theta_s \right) + M_c p- \sigma_0\left(\lambda \right)  
\label{Eq_YieldFunction}
\end{equation}
where, as shown in Fig. \ref{Fig_CohesionTerm},   $M_c$ is the slope of the criterion in the meridional section associated to a Lode's angle $\theta_s=-\frac{\pi}{6}$
and is defined as
\begin{equation}
M_c=\frac{6 \sin \phi}{3-\sin \phi}
\end{equation}
where $\phi$ is the angle of shearing resistance.

\begin{figure}[ht]

\begin{subfigure}{.48\textwidth}
  \centering
  \includegraphics[height=0.9\linewidth]{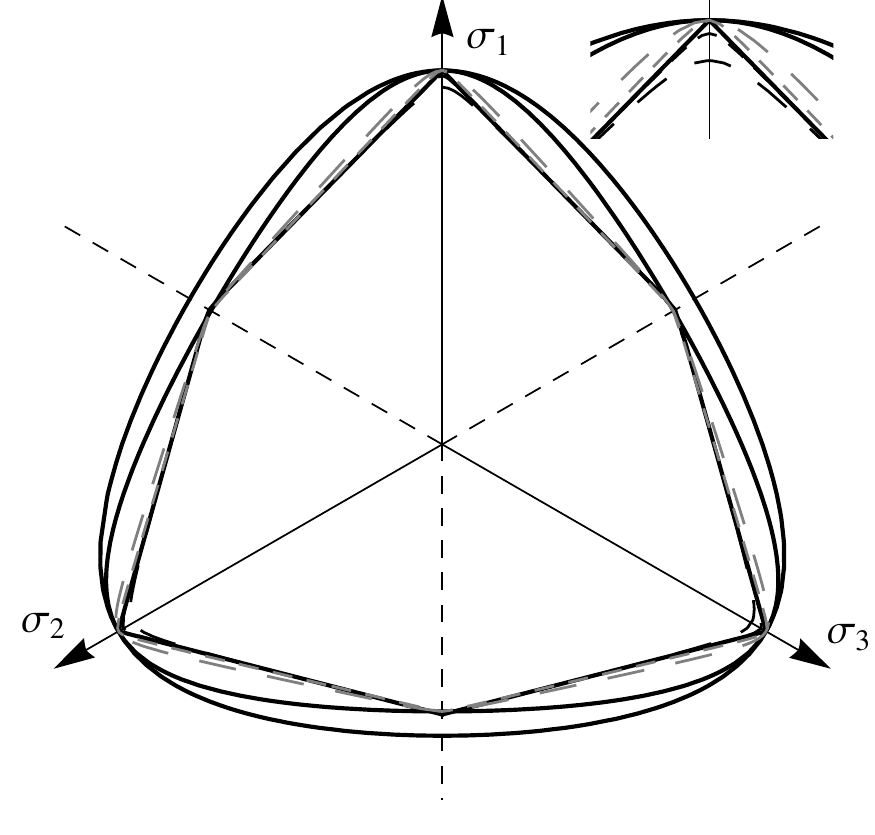}  
  \caption{Mohr-Coulomb,  Matsuoka-Nakai and Lade-Duncan}
  \label{fig:sub-first}
\end{subfigure}
\begin{subfigure}{.48\textwidth}
  \centering
  \includegraphics[height=0.9\linewidth]{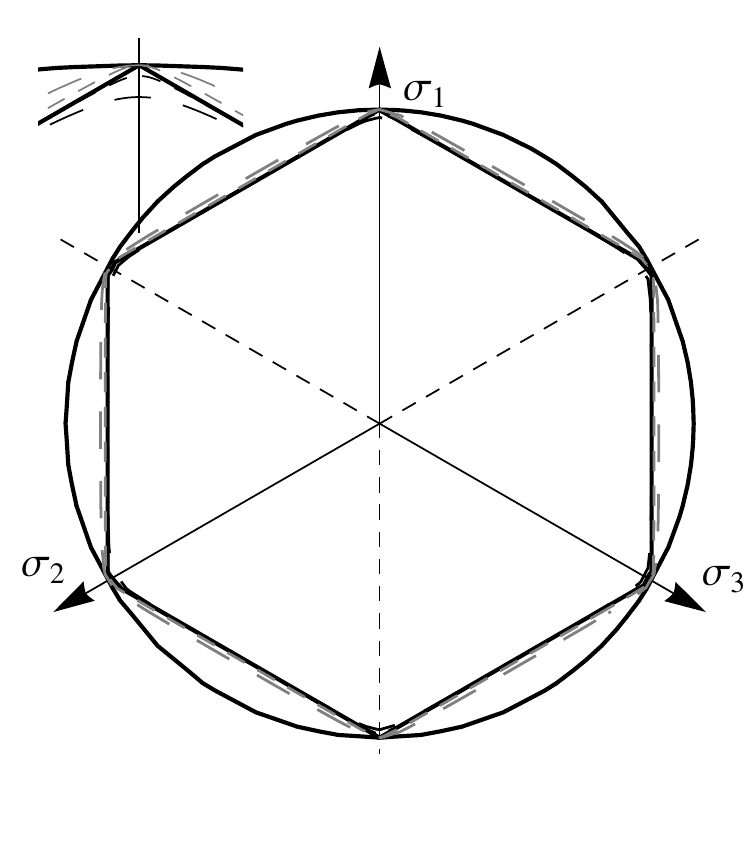}  
  \caption{Tresca,   von Mises and Drucker-Prager}
  \label{fig:sub-second}
\end{subfigure}
\caption{Exemplification of the shapes provided by the $\Gamma(\theta_s)$ function of the GC failure criterion (from Lagioia and Panteghini \cite{LP2016}).  For the Tresca and Mohr-Coulomb criteria rounded versions are also shown.}
\label{Fig_GammaFunctionGraphs}
\end{figure}

The function $\Gamma$ is 
\begin{equation}
\Gamma\left(\theta_s \right)= \alpha_f \, \cos\left[ \frac{\arccos \left(-\beta_f \sin 3 \theta_s \right)}{3}- \gamma_f\frac{\pi}{6}\right]
\label{Eq_GammaFunction}
\end{equation}
and $\alpha_f$, $\beta_f$ and $\gamma_f$ are the three parameters which rule the shape of the yield function in the deviatoric plane and are provided  in Lagioia and Panteghini  \cite{LP2016} without the $f$ subscript (Table \ref{tab:AlfaBetaGamma}).   An additional set of parameters will be also used for the plastic potential surface,  distinguished form the previous one using  a $g$ subscript, so that the possibility is given to adopt surfaces with different deviatoric sections. 
It should also be noted that the minus sign before the $\beta_f$ parameter is introduced in Eq. \ref{Eq_GammaFunction} to account for the sign convention adopted in this paper. 
An instance of the shapes given by the $\Gamma(\theta_s)$ function are shown in Fig. \ref{Fig_GammaFunctionGraphs}.

\begin{figure}[tb]
\centering
\includegraphics[width=0.5\textwidth]{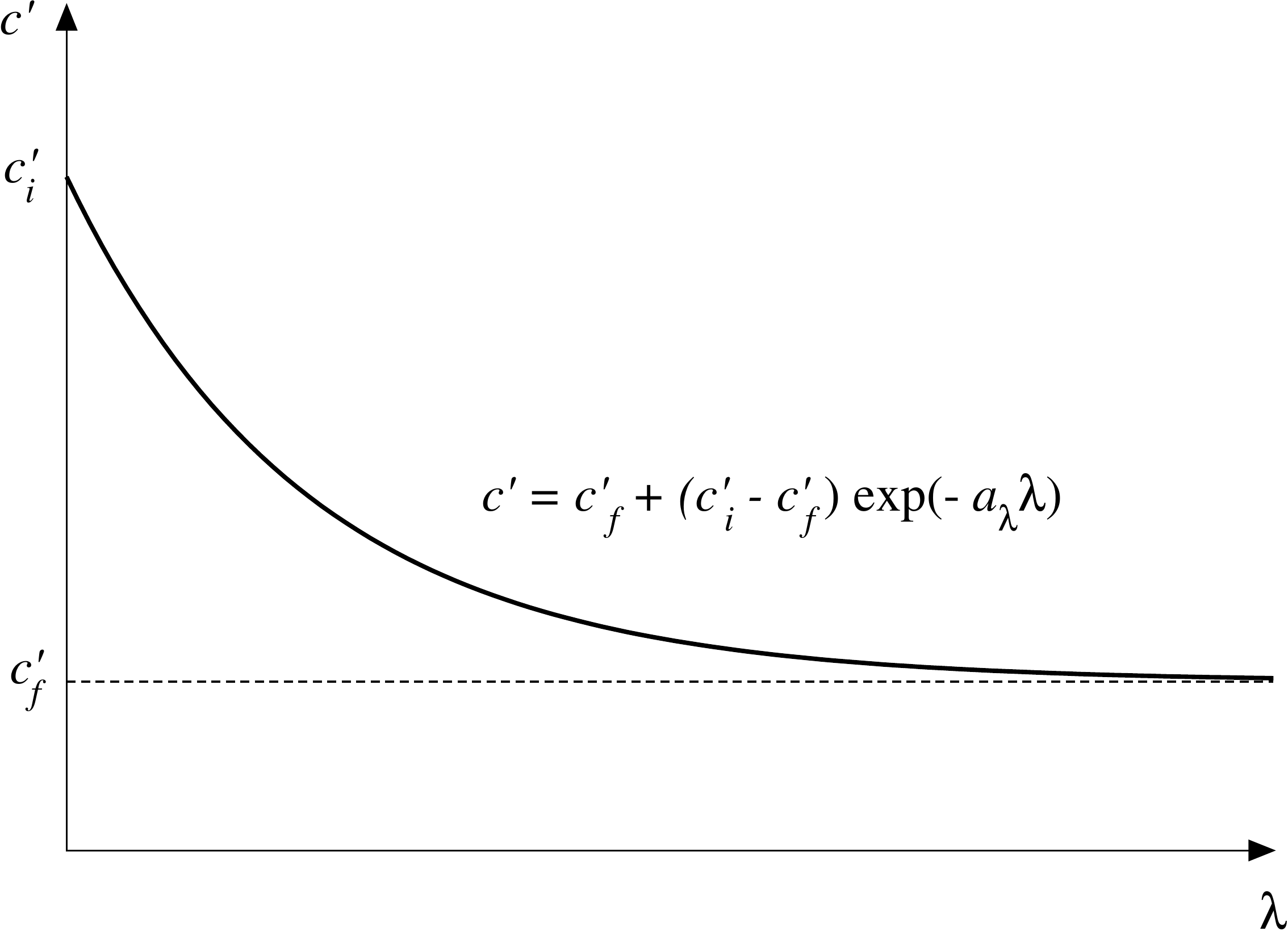}
\caption{Adopted softening of the \textit{cohesion} term of the failure criterion.}
\label{Fig_CohesionSoftening}
\end{figure}

\begin{table*}[!ht]
\begin{tabular}{p{0.15\linewidth}p{0.2\linewidth}p{0.1\linewidth}p{0.45\linewidth}}
\hline
      Model & $\alpha$ & $\beta$ & $\gamma$ \\
    \hline

    von Mises           & 1 & 0 &   1 \\[1. em] 
     Drucker-Prager           & 1 & 0 &   1 \\[1. em] 
     Tresca	             & $\sec\left[\myPiS\right]$ &    1   &  1 \\[1. em] 
   
    \captionfont Mohr-Coulomb        &\captionfont $\sec\left[(\bar{\gamma}+1)\myPiS\right]$ &\captionfont 1    &\captionfont $(1-\bar{\gamma})$ \\[1. em] 
    
    \captionfont Matsuoka-Nakai   &\captionfont $\frac{2}{3}\sqrt{A_1} M_c$ & $\captionfont \frac{A_2}{A_1^{3/2}}$  & \captionfont 0 \\[1. em] 
    \captionfont Lade-Duncan  &  $\captionfont \frac{2}{3}\sqrt{A_1} M_c$ & $\captionfont \frac{A_2}{A_1^{3/2}}$ & \captionfont 0 \\[1. em] 
 
    \hline
    \\[1 em] 
\end{tabular}
\caption{Parameters for the shape function of the generalized criterion of \cite{LP2016}. In addition  $\bar{\gamma}=\frac{6}{\pi} \arctan \frac{\sin \phi}{\sqrt{3}}$,  $\phi$ is the angle of shearing resistance. For Matsuoka-Nakai $A_1=\frac{K_{MN}-3}{K_{MN}-9}$, $A_2=\frac{K_{MN}}{K_{MN}-9}$ and $K_{MN}=\frac{9-\sin^2 \phi}{1-\sin^2 \phi}$, whilst for Lade-Duncan $A_1=\frac{K_{LD}}{K_{LD}-27}$, $A_2=A_1$ and and $K_{LD}=\frac{(3-\sin \phi) ^3}{(1+\sin \phi) ( 1- \sin \phi)^2}$ (modified from \cite{LP2016}). }	
	\label{tab:AlfaBetaGamma}

\end{table*}

The term $\sigma_0$, also shown in Fig. \ref{Fig_CohesionTerm},  accounts for the so called \textit{cohesion} in the classical criterion and is  expressed as 
\begin{equation}
\sigma_0(\lambda)= c'(\lambda) \frac{6 \cos \phi}{3-\sin \phi}
\label{Eq_HardeningParameter}
\end{equation}
It  defines the intercept of the yield/failure criterion in the triaxial compression section (i.e.  $\theta_s=-\frac{\pi}{6}$) with the $q$ axis.
It should be noted that in Eq. \eqref{Eq_HardeningParameter}  $\lambda$ is a generic scalar function of  plastic strains, so that isotropic strain hardening/softening can be accounted for.
In this paper a simple hardening rule will be adopted, defined as
\begin{equation}
c'(\lambda)= c'_f+ (c'_i-c'_f) \, \text{exp} ( -a \lambda) 
\end{equation}
where $c'_i$ and $c'_f$ are the initial and the final cohesions and  $a$ is a parameter which controls the rate of the reduction (Fig. \ref{Fig_CohesionSoftening}).
A similar softening rule is often adopted in numerical analyses of geotechnical problems both in drained and undrained conditions. As an example Summersgill et al.  \cite{Summersgill2018} use a linear reduction of the cohesion term.
In this paper a very non-linear hardening law has been used in order to test the capability of the Cosserat continuum to regularize the boundary value problem even in extreme conditions. It should be noted, however that adopting a more suitable hardening/softening rule can be done without modifying the general structure of the model.

In what follows the plastic potential function will be indicated as 
\begin{equation}
g(p, q, \theta_s)= q \hat{\Gamma}\,(\theta_s) + \hat{M_c} p
\end{equation}
where the $\hat{(\cdot)}$  symbol is introduced to differentiate from the yield function.

Plastic strain tensors are obtained through the plastic potential as 
\begin{equation}
\begin{aligned}
\dot{\tenss{\gamma}}^p= \der{g(p,q,\theta_s)}{\tenss{\sigma}} \dot{ \lambda} \\
\dot{\tenss{\chi}}^p= \der{g(p,q,\theta_s)}{\tenss{\mu}} \dot{ \lambda}
\end{aligned}
\end{equation}
where $\dot \lambda$ is the plastic multiplier and the standard  assumption is made that the strain and curvature increments can be split into their elastic and plastic components
\begin{equation}
\begin{aligned}
\dot{\tenss{\gamma}}= \dot{\tenss{\gamma}}^e+\dot{\tenss{\gamma}}^p\\
\dot{\tenss{\chi}}= \dot{\tenss{\chi}}^e+\dot{\tenss{\chi}}^p
\end{aligned}
\end{equation}

\section{Conclusions}
An elasto-plastic constitutive model for the linear formulation of the  Cosserat continuum has been presented.
The model features non-associated flow and hardening/softening behaviour, whilst linear hyper-elasticity is adopted to reproduce the recoverable response.

For the formulation  of the yield and plastic potential functions,  the definition of the  \textit{equivalent von Mises stress}  used in the Cauchy continuum has been extended to the Cosserat material exploiting   Hencky's  \cite{Hencky1924} interpretation of the von Mises criterion.  The resulting expression for the equivalent von Mises stress coincides with that of Muehlhaus and Vardoulakis  \cite{Muehlhaus1987} which was obtained on the basis of  considerations of micro-mechanics  and averaging. 
The same expression can also be obtained using the theory of representations.

The key feature of the model is the capability of providing different shapes of the yield and plastic potential surfaces in the deviatoric plane.  Whilst most models in the literature adopt surfaces with a circular shape in the deviatoric plane, hence dropping the dependency on the Lode's angle,  a general shape function has been used in the proposed model for the deviatoric section.
The dependency of the those functions on the Lode's angle has hence been introduced in the constitutive model, which makes it significant for practical applications.
Both surfaces are defined using the function,  proposed by Lagioia and Panteghini \cite{LP2016},  which is an exact definition of most classical yield and failure criteria.

\bibliographystyle{plain}
\bibliography{short,paper}

\end{document}